\documentclass[11pt,reqno]{amsart}
\usepackage{hyperref}
\usepackage{graphicx}
\usepackage{color}
\usepackage{amsfonts,amssymb}
\usepackage{bbm} 
\usepackage[latin9]{inputenc}

\usepackage{mathrsfs}

\usepackage{a4wide}

\usepackage{color} 

\definecolor{green}{rgb}{0,0.6,0}
\definecolor{blue}{rgb}{0,0,1}

\hypersetup{colorlinks,
linkcolor=blue,
filecolor=green,
citecolor=green}

\theoremstyle{plain}

\newtheorem{neu}{}[section]
\newtheorem{Cor}[neu]{Corollary}
\newtheorem*{Cor*}{Corollary}
\newtheorem{Thm}[neu]{Theorem}
\newtheorem*{Thm*}{Theorem}
\newtheorem{Prop}[neu]{Proposition}
\newtheorem*{Prop*}{Proposition}
\theoremstyle{definition}
\newtheorem{Lemma}[neu]{Lemma}
\newtheorem*{Rmk*}{Remark}
\newtheorem{Rmk}[neu]{Remark}

\newtheorem*{Ex*}{Example}
\newtheorem{Qu}[neu]{Question}
\newtheorem*{Qu*}{Question}

\newtheorem{Def}[neu]{Definition}

\theoremstyle{remark}

\theoremstyle{definition}

\newcommand{\p}{\partial}
\newcommand{\om}{\omega}

\newcommand{\pf}{\longrightarrow}

\newcommand{\N}{{\mathbb{N}}}
\newcommand{\Z}{{\mathbb{Z}}}
\newcommand{\R}{{\mathbb{R}}}

\newcommand{\M}{\mathcal{M}}
\newcommand{\LLL}{\mathscr{L}}

\newcommand{\B}{{\bf B}}

\renewcommand{\H}{{\bf H}}

\renewcommand{\H}{\mathrm{H}}

\newcommand{\CF}{\mathrm{CF}}
\newcommand{\HF}{\mathrm{HF}}
\newcommand{\RFH}{\mathrm{RFH}}
\newcommand{\CM}{\mathrm{CM}}

\newcommand{\ind}{{\rm ind}}

\newcommand{\Crit}{{\rm Crit}}

\newcommand{\codim}{{\rm codim}}

\newcommand{\Ham}{\mathrm{Ham}}
\newcommand{\Supp}{\mathrm{Supp}}

\newcommand{\EE}{\mathcal{E}}
\newcommand{\BB}{\mathcal{B}}

\newcommand{\MM}{\mathcal{M}}

\newcommand{\LL}{\mathcal{L}}
\newcommand{\FF}{\mathcal{F}}

\renewcommand{\AA}{\mathcal{A}}

\newcommand{\comment}[1]{}
\newcommand{\x}{\times}

\newcommand{\beq}{\begin{equation}}
\newcommand{\beqn}{\begin{equation}\nonumber}
\newcommand{\eeq}{\end{equation}}

\newcommand{\bea}{\begin{equation}\begin{aligned}}
\newcommand{\bean}{\begin{equation}\begin{aligned}\nonumber}
\newcommand{\eea}{\end{aligned}\end{equation}}

\numberwithin{equation}{section}
\begin{document}
\title{K\"unneth formula in Rabinowitz Floer homology}
\author{Jungsoo Kang}

\address{Department of Mathematics, Seoul National University,
Kwanakgu Shinrim, San56-1 Seoul, South Korea, Email:
hoho159@snu.ac.kr}

\begin{abstract}
Rabinowitz Floer homology has been investigated on submanifolds of contact type. The contact condition, however, is quite restrictive. For example, a
product of contact hypersurfaces is rarely of contact type. In this article, we study Rabinowitz Floer homology for product manifolds which are not necessarily of contact type. We show for a class of product manifolds that there are infinitely many leafwise intersection points by proving the K\"unneth formula for Rabinowitz Floer homology.
\end{abstract}
\maketitle

\section{Introduction}
Rabinowitz Floer homology has been extensively studied in recent years because of its interrelation with the leafwise intersection problem. However Rabinowitz Floer homology (to be honest, the perturbed Rabinowitz action functional) has worked principally on a contact submanifold and little research has been conducted on a non-contact case. Our primary objective in this paper is to find leafwise intersection points and define Rabinowitz Floer homology for this class of submanifolds which are not necessarily of contact type. In addition we show for the class that there are infinitely many leafwise intersection points by proving the K\"unneth formula for Rabinowitz Floer homology. For simplicity, throughout this paper, we use $\Z/2$-coefficients for Rabinowitz Floer homology, but expect the K\"unneth formula continues to hold with $\Z$-coefficient.\\[-1.5ex]

We consider restricted contact hypersurfaces $(\Sigma_1,\lambda_1)$ resp. $(\Sigma_2,\lambda_2)$ in exact symplectic manifolds $(M_1,\om_1=d\lambda_1)$ resp. $(M_2,\om_2=d\lambda_2)$. Moreover we assume that $\Sigma_1$ resp. $\Sigma_2$ bounds a compact region in $M_1$ resp. $M_2$ and that those $M_1$ and $M_2$ are convex at infinity; that is, they are symplectomorphic to the symplectization of a compact contact manifold at infinity. Given $F_1\in(S^1\x M_1)$, $F_2\in C^\infty(S^1\x M_2)$, the operation
\beqn
\big(F_1\oplus F_2\big)(t,x,y)=F_1(t,x)+F_2(t,y), \qquad(t,x,y)\in S^1\x M_1\x M_2
\eeq
provides a time-dependent Hamiltonian function $F_1\oplus F_2\in C^\infty(S^1\x M_1\x M_2)$. We also introduce projection maps $\pi_1:M_1\x M_2\to M_1$ and $\pi_2:M_1\x M_2\to M_2$; then $(M_1\x M_2,\om_1\oplus\om_2)$ admits the symplectic structure $\om_1\oplus\om_2=\pi_1^*\om_1+\pi_2^*\om_2$.\\[-1.5ex]

On $(\Sigma_1\x\Sigma_2,M_1\x M_2)$, we define the perturbed Rabinowitz action functional $\AA^{\widetilde H_1,\widetilde H_2}_F$ as in  \eqref{eq:Rabinowitz action functional}. Since $\Sigma_1\x \Sigma_2$ is a stable submanifold, we can define Floer homology of $\AA^{\widetilde H_1,\widetilde H_2}$ when $F\equiv0$ (refer to \cite{Ka2} for definitions and constructions). This Floer homology $\HF(\AA^{\widetilde H_1,\widetilde H_2})$ is called {\em Rabinowitz Floer homology} and denoted by $\RFH(\Sigma_1\x\Sigma_2,M_1\x M_2)$, see Section 3. By the standard continuation method in Floer theory, $\HF(\AA^{\widetilde H_1,\widetilde H_2}_F)$ and $\RFH(\Sigma_1\x\Sigma_2,M_1\x M_2)$ are isomorphic whenever $\HF(\AA^{\widetilde H_1,\widetilde H_2}_F)$ is defined.

\noindent\textbf{Theorem A.}
{\em The Floer homologies $\RFH(\Sigma_1\x\Sigma_2,M_1\x M_2)$ and $\HF(\AA^{\widetilde H_1,\widetilde H_2}_{F_1\oplus F_2})$ are well-defined. Moreover, we have the following K\"unneth formula in Rabinowitz Floer homology:}
\beqn
\RFH_n(\Sigma_1\x\Sigma_2,M_1\x M_2)\cong\bigoplus_{p=0}^n\RFH_p(\Sigma_1,M_1)\otimes\RFH_{n-p}(\Sigma_2,M_2).\\[1ex]
\eeq
Here, $\RFH_p(\Sigma_1,M_1)$ (resp. $\RFH_{n-p}(\Sigma_2,M_2)$) is the  Rabinowitz Floer homology for the restricted contact hypersurface $\Sigma_1$ in $M_1$ (resp. $\Sigma_2$ in $M_2$), see \cite{AF1} or Section 3.

\begin{Rmk}
In this paper, we unfortunately establish compactness of gradient flow lines of the Rabinowitz action functional only for perturbations of the form $F=F_1\oplus F_2$. Thus we cannot study the existence problem of leafwise intersection points for an arbitrary perturbation. However, if $\Sigma_1\x\Sigma_2$ has contact type in the sense of Bolle \cite{Bo1,Bo2} (see Section 4), the Floer homology $\HF(\AA^{\widetilde H_1,\widetilde H_2}_F)$ is defined for all perturbations, see \cite{Ka2}. We note that, in general, $\Sigma_1\x\Sigma_2$ is not of contact type in the sense of Bolle. For example, $S^3\x S^3$ is not a contact submanifold in $\R^8$, see Remark \ref{rmk:restriction of contact}.
\end{Rmk}

\begin{Qu}
What perturbations have a leafwise intersection point on $(\Sigma_1\x\Sigma_2,M_1\x M_2)$?
\end{Qu}

\begin{Rmk}
Once one verifies compactness of gradient flow lines of the Rabinowitz action functional for a given perturbation $F$, it guarantees the existence of leafwise intersection points for that $F$ by using the stretching the neck argument in \cite{AF1}. In this paper, we are able to compactify gradient flow lines of $\AA^{\widetilde H_1,\widetilde H_2}_{F_1\oplus F_2}$, and thus it guarantees the existence of leafwise intersection points of $F_1\oplus F_2$; but, this directly follows from the result in \cite{AF1} that each $F_1$ and $F_2$ has a leafwise intersection point on $\Sigma_1$ and $\Sigma_2$ respectively.
\end{Rmk}
\begin{Def}
The {\em Hamiltonian vector field} $X_F$ on a symplectic manifold $(M,\om)$ is defined explicitly by $i_{X_F}\om=dF$ for a Hamiltonian function $F\in C^\infty(S^1\x M)$, and we call its time one flow $\phi_F$ the {\em Hamiltonian diffeomorphism}. We denote by $\Ham_c(M,\om)$ the group of Hamiltonian diffeomorphisms generated by compactly supported Hamiltonian function. This group has a well-known norm introduced by Hofer (see Definition \ref{def:Ham}).
\end{Def}

\begin{Def}
We denote by $\wp(\Sigma_1,\lambda_1)>0$ the {\em minimal period} of closed Reeb orbits of $(\Sigma_1,\lambda_1)$ which are contractible in $M_1$. If there is no contractible closed Reeb orbit we set $\wp(\Sigma_1,\lambda_1)=\infty$.
\end{Def}
In Theorem B we do not consider $\Sigma_2$, and $M_2$ need to be closed.\\[-1ex]

\noindent\textbf{Theorem B.} {\em Let $(M_2,\om_2)$ be a closed and symplectically aspherical, i.e. $\om_2|_{\pi_2(M_2)}=0$, symplectic manifold. Then, although $\Sigma_1\x M_2$ is not a contact hypersurface,}
\begin{enumerate}
\item[(B1)] {\em $\Sigma_1\x M_2$ has a leafwise intersection point for $\phi\in\Ham_c(M_1\x M_2,\om_1\oplus\om_2)$ with Hofer-norm $||\phi||<\wp(\Sigma_1,\lambda_1)$ even if $\Sigma_1$ does not bound a compact region in $M_1$.
\item[(B2)] The Rabinowitz Floer homology $\RFH(\Sigma_1\x M_2,M_1\x M_2)$ can be defined when $\Sigma_1$ bounds a compact region in $M_1$. Moreover, we have the K\"unneth formula:}
\beqn
\RFH_n(\Sigma_1\x M_2,M_1\x M_2)\cong\bigoplus_{p=0}^n\RFH_p(\Sigma_1,M_1)\otimes\H_{n-p}(M_2).\\[1ex]
\eeq
\end{enumerate}

To prove Theorem B without any contact conditions, we need to show a special version of isoperimetric inequality, see Lemma \eqref{lemma:isoperimetric ineq}.
\begin{Rmk}
It is worth emphasizing that $\Sigma_1\x M_2$ is not necessarily of restricted contact type. For instance, if $M_2$ is not exact, then $\Sigma_1\x M_2$ is never of restricted contact type. Nevertheless, interestingly enough, we can achieve compactness of gradient flow lines of Rabinowitz action functional for an arbitrary perturbation $F\in\Ham_c(M_1\x M_2,\om_1\oplus\om_2)$; accordingly Floer homology of the Rabinowitz action functional with any perturbations is well-defined.
\end{Rmk}

The K\"unneth formula enable us to compute the Rabinowitz Floer homology of a product manifold in terms of Rabinowitz Floer homology of each manifolds. As applications, in Section 4 we shall prove the following two corollaries.\\[1ex]
\noindent\textbf{Corollary A.}
{\em Let $N$ be a closed Riemannian manifold of $\dim N\geq2$ with $\dim\H_*(\Lambda N)=\infty$ where $\Lambda N$ is the free loop space of $N$. Then there exists infinitely many leafwise intersection points for a generic $\phi\in\Ham_c(T^*S^1\x T^*N)$ on $(S^*S^1\x S^*N,T^*S^1\x T^*N)$.}

\begin{Rmk}
Since $(S^*S^1\x S^*N,T^*S^1\x T^*N)$ is of restricted contact type in the sense of Bolle (Lemma \ref{lemma:SS1 SN is contact}), $\phi$ in Corollary A is not necessarily of product type. If $\phi$ has product type, then the above result is obvious by \cite{AF1,AF2}. Unlike Corollary A, the following Corollary B does not assume the contact condition since Theorem B does not need any contact conditions.
\end{Rmk}

\noindent\textbf{Corollary B.}
{\em Let $M$ be a closed and symplectically aspherical symplectic manifold and $N$ be as above. Then a generic $\phi\in\Ham_c(T^*N\x M)$ has infinitely many leafwise intersection points on $(S^*N\x M,T^*N\x M)$.}

\begin{Rmk}
If $\pi_1(N)$ is finite then $\dim\H_*(\Lambda N)=\infty$ by \cite{Sullivan_Vigue_the_homology_theory_of_the_closed_geodesic_problem}. If the number of conjugacy classes of $\pi_1(N)$ is infinite then $\dim\H_0(\Lambda N)=\infty$. Therefore, the only remaining case is if $\pi_1(N)$ is infinite but the number of conjugacy classes of $\pi_1(N)$ is finite.
\end{Rmk}
\subsection{Leafwise intersections}
Let $(M,\om)$ be a $2n$ dimensional symplectic manifold and $\Sigma$ be a coisotropic submanifold of codimension $0\leq k\leq n$. Then the symplectic structure $\om$ determines a symplectic orthogonal bundle $T\Sigma^\om\subset T\Sigma$ as follows:
\beqn
T\Sigma^\om:=\{(x,\xi)\in T\Sigma\,|\,\om_x(\xi,\zeta)=0 \textrm{ for all } \zeta\in T_x\Sigma\}
\eeq
Since $\om$ is closed, $T\Sigma^\om$ is integrable, thus $\Sigma$ is foliated by the leaves of the characteristic foliation and we denote by $L_x$ the isotropic leaf through $x$. We call $x\in\Sigma$ a {\em leafwise intersection point} of $\phi\in\Ham(M,\om)$ if $x\in L_x\cap\phi(L_x)$. In the extremal case $k=n$, Lagrangian submanifold consists of only one leaf. Thus a leafwise intersection point is nothing but a Lagrangian intersection point in the case $k=n$. In the other extremal case that $k=0$, a leafwise intersection corresponds to a periodic orbit of $\phi$.\\[-1ex]

The leafwise intersection problem was initiated by Moser \cite{M} and pursued further in \cite{Ba,H,EH,Gi1,Dr,Gu,AF1,AF2,AF3,AF4,Z,Ka1,Ka2,Me,AMc,AM}. We refer to \cite{AF1,Ka2} for the history of the problem. In particular, Albers-Frauenfelder approached the problem by means of the perturbed Rabinowitz action functional and much relevant research has been conducted in \cite{AF1,AF2,AF3,AF4,CF,CFO,CFP,AS,Ka1,Ka2,Me,AM}. We refer to \cite{AF5} for a brief survey on Rabinowitz Floer theory.

\section{Rabinowitz action functional on product manifolds}

Since $\Sigma_1$ and $\Sigma_2$ are contact hypersurfaces, there exist associated Liouville vector fields $Y_1$ resp. $Y_2$ on $M_1$ resp. $M_2$ such that $\LL_{Y_i}\om_i=\om_i$ and $Y_i\pitchfork\Sigma_i$ for $i=1,2$. We denote by $\phi_{Y_i}^t$ the flow of $Y_i$ and fix $\delta>0$ such that $\phi_{Y_i}^t|_{\Sigma_i}$ is defined for $|t|<\delta$. Since $\Sigma_1$ resp. $\Sigma_2$ bounds a compact region in $M_1$ resp. $M_2$, we are able to define Hamiltonian functions $G_1\in C^\infty(M_1)$ and $G_2\in C^\infty(M_2)$ so that
\begin{enumerate}
\item $G_1^{-1}(0)=\Sigma_1$ and $G_2^{-1}(0)=\Sigma_2$ are regular level sets;
\item $dG_1$ and $dG_2$ have compact supports;
\item $G_i(\phi_{Y_i}^t(x_i))=t$ for all $x_i\in\Sigma_i$, $i=1,2$, and $|t|<\delta$;
\end{enumerate}
We extend $G_1$, $G_2$ to be defined on the whole of $M_1\x M_2$:
\bean
\widetilde G_i:M_1\x M_2&\pf\R\qquad i=1,2\\
(x_1,x_2)&\longmapsto G_i(x_i).
\eea

\begin{Def}\label{def:Moser triple}
Given time-dependent Hamiltonian functions $\widetilde H_1,\widetilde H_2,F\in C^\infty(S^1\x M_1\x M_2)$, a triple $(\widetilde H_1,\widetilde H_2,F)$ is called a {\em Moser triple} if it satisfies
\begin{enumerate}
\item their time supports are disjoint, i.e.
\beqn
\widetilde H_1(t,\cdot)=\widetilde H_2(t,\cdot)=0 \quad\textrm{for}\,\,\, \forall t\in[0,\frac{1}{2}] \quad\textrm{and}\quad F(t,\cdot)=0 \quad\textrm{for}\,\,\, \forall t\in[\frac{1}{2},1].
\eeq
\item $F=F_1\oplus F_2$ for some $F_1\in C_c^\infty(S^1\x M_1)$, $F_2\in C_c^\infty(S^1\x M_2)$.
\item $\widetilde H_1$ and $\widetilde H_2$ are weakly time-dependent Hamiltonian functions. That is, $\widetilde H_1$ and $\widetilde H_2$ are of the form $\big(\widetilde H_1(t,x),\widetilde H_2(t,x)\big)=\chi(t)\big(\widetilde G_1(x),\widetilde G_2(x)\big)$ for $\chi:S^1\to S^1$ with $\int_0^1\chi dt=1$ and $\Supp\chi\subset(\frac{1}{2},1)$.
\end{enumerate}
\end{Def}

Next, we recall the definition of the Hofer norm.
\begin{Def}\label{def:Ham}
Let $F\in C^\infty_c(S^1\x M,\mathbb{R})$ be a compactly supported time-dependent Hamiltonian function on a symplectic manifold $(M,\om)$. We set
\beqn
||F||_+:=\int_0^1\max_{x\in M} F(t,x) dt\qquad||F||_-:=-\int_0^1\min_{x\in M} F(t,x) dt=||-F||_+
\eeq
and
\beqn
||F||=||F||_++||F||_-.\;
\eeq
For $\phi\in \Ham_c(M,\om)$ the Hofer norm is
\beqn
||\phi||=\inf\{||F||\mid \phi=\phi_F, F\in C^\infty_c(S^1\x M,\mathbb{R})\}.\;
\eeq

\end{Def}

\begin{Lemma}\label{lemma:norms equivalent}
For all $\phi\in\Ham_c(M,\om)$
\beqn
||\phi||=|||\phi|||:=\inf\big\{||F||\mid \phi=\phi_F,\;F(t,\cdot)=0\;\;\forall t\in[\tfrac12,1]\big\}\;.
\eeq
\end{Lemma}
\begin{proof}
To prove $||\phi||\geq|||\phi|||$, pick a smooth monotone increasing map $r:[0,1]\to[0,1]$ with $r(0)=0$ and $r(1/2)=1$. For $F$ with $\phi_F=\phi$ we set $F^r(t,x):=r'(t)F(r(t),x)$. Then a direct computation shows $\phi_{F^r}=\phi_F$, $||F^r||=||F||$, and $F^r(t,x)=0$ for all $t\in[\tfrac12,1]$. The reverse inequality is obvious.
\end{proof}

We denote by $\LLL=\LLL_{M_1\x M_2}\subset C^\infty(S^1,M_1\x M_2)$ the component of contractible loops in $M_1\x M_2$. With a Moser triple $(H_1,H_2,F)$, the perturbed Rabinowitz action functional $\AA_{F}^{\widetilde H_1,\widetilde H_2}(v,\eta_1,\eta_2):\LLL\times\R^2\longrightarrow\mathbb{R}$ is defined as follows:

\beq\label{eq:Rabinowitz action functional}
\AA_{F}^{\widetilde H_1,\widetilde H_2}(v,\eta_1,\eta_2)=-\int_0^1 v^*\lambda_1\oplus\lambda_2-\eta_1\int_0^1\widetilde H_1(t,v)dt-\eta_2\int_0^1\widetilde H_2(t,v)dt-\int_0^1F(t,v)dt
\eeq
where $\lambda_1\oplus\lambda_2=\pi_1^*\lambda_1+\pi_2^*\lambda_2$. The real numbers $\eta_1$ and $\eta_2$ can be thought of as Lagrange multipliers.

Critical points $(v,\eta_1,\eta_2)\in\Crit\AA_{F}^{\widetilde H_1,\widetilde H_2}$ satisfy
\beq\label{eqn:critical point equation}\left.
\begin{aligned}
&\partial_tv=\eta_1X_{\widetilde H_1}(t,v)+\eta_2X_{\widetilde H_2}(t,v)+X_{F}(t,v),\\[1ex]
&\int_0^1\widetilde H_1(t,v)dt=0,\\[1ex]
&\int_0^1\widetilde H_2(t,v)dt=0.
\end{aligned}
\;\;\right\}
\eeq

Albers-Frauenfelder \cite{AF1} observed that a critical point of $\AA_{F}^{\widetilde H_1,\widetilde H_2}$ gives rise to a leafwise intersection point. (In fact, they proved the following proposition for the codimensional one case, yet their proof continues to hold in our case, see \cite{Ka2} also.)

\begin{Def}\label{def:periodic_LI}
A leafwise coisotropic intersection point $x\in\Sigma_1\x\Sigma_2$ is called {\em periodic} if the leaf $L_x$ contains neither a closed Reeb orbit in $\Sigma_1$ nor a closed Reeb orbit in $\Sigma_2$.
\end{Def}

\begin{Prop}\label{prop:critical point answers question}
\cite{AF1} Let $(v,\eta_1,\eta_2)\in\Crit\AA_{F}^{\widetilde H_1,\widetilde H_2}$. Then $x=v(1/2)$ satisfies $\phi_F(x)\in L_x$. Thus, $x$ is a leafwise intersection point. Moreover, the map
\beqn
\Crit\AA^{\widetilde H_1,\widetilde H_2}_F\pf\big\{\textrm{leafwise intersections}\big\}
\eeq
is injective unless there exists a periodic leafwise intersection.
\end{Prop}

We choose a compatible almost complex structure $J_1$ on $M_1$ and define the metric on $(M_1,\om_1)$ by $g_1(\cdot,\cdot)=\om_1(\cdot,J_1\cdot)$. Analogously we also define the metric on $(M_2,\om_2)$, $g_2(\cdot,\cdot)=\om_2(\cdot,J_2\cdot)$. Then $g=g_1\oplus g_2$ which is the metric on $(M_1\x M_2,\om_1\oplus\om_2)$ induces a metric $m$ on the tangent space $T_{(v,\eta_1,\eta_2)}(\LLL\times\R^2)\cong T_v\LLL\x\R^2$ as follows:
\beqn
m_{(v,\eta_1,\eta_2)}\big((\hat v^1,\hat\eta^1_1,\hat\eta_2^1),(\hat v^2,\hat\eta_1^2,\hat\eta^2_2)\big):=\int_0^1g_v(\hat v^1,\hat v^2)dt+\hat\eta^1_1\hat\eta^2_1+\hat\eta^1_2\hat\eta^2_2\;.
\eeq

\begin{Def}
A map $w=(v,\eta_1,\eta_2)\in C^\infty(\mathbb{R},\LLL\times\mathbb{R}^2)$ which solves
\beq\label{eqn:gradient flow line}
\partial_s w(s)+\nabla_m \AA_{F}^{\widetilde H_1,\widetilde H_2}(w(s))=0\;
\eeq
is called a gradient flow line of $\AA_{F}^{\widetilde H_1,\widetilde H_2}$ with respect to the metric $m$.
\end{Def}

According to Floer's interpretation, the gradient flow equation \eqref{eqn:gradient flow line} can be interpreted as maps $v(s,t):\mathbb{R}\x S^1\to M_1\x M_2$ and $\eta_1(s)$, $\eta_2(s):\R\to\R$ solving

\beq\label{eqn:gradient flow equation}\left.
\begin{aligned}
&\partial_sv+J(v)\bigr(\partial_tv-\eta_1X_{\widetilde H_1}(t,v)-\eta_2X_{\widetilde H_2}(t,v)-X_{F}(t,v)\bigr)=0,\\[1ex]
&\partial_s\eta_1-\int_0^1\widetilde H_1(t,v)dt=0,\\[1ex]
&\partial_s\eta_2-\int_0^1\widetilde H_2(t,v)dt=0 .
\end{aligned}
\;\;\right\}
\eeq

\begin{Def}
The energy of a map $w\in C^{\infty}(\R,\LLL\times\R^2)$ is defined by
\beqn
E(w):=\int_{-\infty}^\infty||\partial_s w||_m^2ds\;.
\eeq
\end{Def}

\begin{Lemma}\label{lemma:energy estimate for gradient lines}
Let $w$ be a gradient flow line of $\AA_{F}^{\widetilde H_1,\widetilde H_2}$. Then
\beq\label{eqn:energy estimate for gradient lines}
E(w)=\AA_{F}^{\widetilde H_1,\widetilde H_2}(w_-)-\AA_{F}^{H_1,H_2}(w_+)\;.
\eeq
where $w_\pm=\lim_{s\to\pm\infty}w(s)$.
\end{Lemma}
\begin{proof}
It follows from the gradient flow equation \eqref{eqn:gradient flow line}.
\bean
E(w)&=\int_{-\infty}^{\infty}m\bigr(-\nabla_m\AA^{\widetilde H_1,\widetilde H_2}_{F}(w(s)),\p_s w(s)\bigr)ds\\
&=-\int_{-\infty}^\infty d\AA_{F}^{\widetilde H_1,\widetilde H_2}(w(s))(\partial_sw(s))ds\\[1ex]
&=-\int_{-\infty}^\infty \frac{d}{ds}\Big(\AA_{F}^{\widetilde H_1,\widetilde H_2}(w(s))\Big)ds\\[1ex]
&=\AA_{F}^{\widetilde H_1,\widetilde H_2}(w_-)-\AA_{F}^{\widetilde H_1,\widetilde H_2}(w_+).
\eea
\end{proof}

\subsection{Compactness of gradient flow lines}
In order to define Rabinowitz Floer homology, we need compactness of gradient flow lines of the Rabinowitz action functional with fixed asymptotic data. More specifically we show the following theorem. In the rest of this section, our perturbation $F\in C^\infty_c(S^1\x M_1\x M_2)$ is of the form $F_1\oplus F_2$ for some $F_1\in C^\infty_c(S^1\x M_1)$ and $F_2\in C^\infty_c(S^1\x M_2)$.
\begin{Thm}\label{thm:compactness}
Let $\{w_n\}_{n\in\N}$ be a sequence of gradient flow lines of $\AA^{\widetilde H_1,\widetilde H_2}_F$ for which there exist $a<b$ such that
\beqn
a\leq\AA^{\widetilde H_1,\widetilde H_2}_F(w_n(s))\leq b,\qquad \textrm{ for all } s\in\R.
\eeq
Then for every reparametization sequence $\sigma_n\in\R$, the sequence $w_n(\cdot+\sigma_n)$ has a subsequence which is converges in $C^\infty_\mathrm{loc}(\R,\LLL\x\R^2)$.
\end{Thm}
\begin{proof}
In order to prove the theorem, we need to verify the following three ingredients.
\begin{enumerate}
\item a uniform $L^\infty$-bound on $v_n$,
\item a uniform $L^\infty$-bound on $\eta_{1n}$, $\eta_{2n}$,
\item a uniform $L^\infty$-bound the derivatives of $v_n$.
\end{enumerate}
for a sequence of gradient flow lines $\{(v_n,\eta_{1n},\eta_{2n})\}_{n\in\N}$. Once we establish (2), the others follow by standard arguments in Floer theory. At the end of this section, we prove Theorem \ref{thm:bound on eta} which proves (2) and thus completes the proof of Theorem \ref{thm:compactness}.
\end{proof}

First of all, we introduce two auxiliary action functionals $\AA_1,\AA_2:\LLL_{M_1\x M_2}\x\R^2\pf\R$:
\bean
\AA_1(v,\eta_1,\eta_2):=\int_0^1 v^*\pi_1^*\lambda_1-\eta_1\int_0^1H_1(t,v)dt-\int_0^1F(t,v)dt,\\
\AA_2(v,\eta_1,\eta_2):=\int_0^1 v^*\pi_2^*\lambda_2-\eta_2\int_0^1H_2(t,v)dt-\int_0^1F(t,v)dt.
\eea

\begin{Lemma}\label{lemma:uniform bound of AA_i}
Let $w=(v,\eta_1,\eta_2)\in C^\infty(\R,\LLL\x\R^2)$ be a gradient flow line of $\AA^{\widetilde H_1,\widetilde H_2}_F$ with asymptotic ends $w_-=(v_-,\eta_{1-},\eta_{2-})$ and  $w_+=(v_+,\eta_{1+},\eta_{2+})$. Then the  action values of $\AA_1$ and $\AA_2$ are bounded along $w$ in terms of the asymptotic data:
\bean
\textrm{(i)}\quad\AA_1(w(s))&\leq2|\AA_1(w_-)|+|\AA_1(w_+)|+4||F_2||_+,\quad\forall s\in\R;\\[1ex]
\textrm{(ii)}\quad\AA_2(w(s))&\leq2|\AA_2(w_-)|+|\AA_2(w_+)|+4||F_1||_+,\quad\forall s\in\R.
\eea
\end{Lemma}
\begin{proof}
We only show the first inequality, the later one is proved in a similar way. Since it holds that $\pi_{1*}X_F=X_{F_1}$, $\pi_{2*}X_F=X_{F_2}$, and $i_{X_{\widetilde H_2}}\pi_1^*\om_1=0$, we compute
\bean
\frac{d}{ds}\AA_1(w(s))&=d\AA_1(w(s))[\p_s w(s)]\\
&=\int_0^1\pi_1^*\om_1\bigr(\p_tv,\p_s v)-\int_0^1\om_1\oplus\om_2\big(\eta_1X_{H_1}(t,v)+X_{F}(t,v),\p_s v\big)\\
&\quad-\Big(\int_0^1\widetilde H_1(t,v)dt\Big)^2\\
&=\int_0^1\pi_1^*\om_1\bigr(\p_tv-\eta_1X_{\widetilde H_1}(t,v)-X_{F}(t,v),\p_s v\bigr)dt\\
&\quad-\int_0^1\pi_2^*\om_2(X_{F}(t,v),\p_s v)dt-\Big(\int_0^1 \widetilde H_1(t,v)dt\Big)^2\\
&=-\int_0^1\pi_1^*\om_1(\p_sv,J\p_sv)dt-\int_0^1\frac{d}{ds}F_2(t,\pi_2\circ v)dt-\Big(\int_0^1\widetilde H_1(t,v)dt\Big)^2.
\eea
Integrating the above equality from $-\infty$ to any $s_0\in\R$, we have
\bea\label{eq:AA_1}
\AA_1(w(s_0))-\AA_1(w_-)&=\int_{-\infty}^{s_0}\frac{d}{ds}\AA_1(w(s))ds\\
&=-\int_{-\infty}^{s_0}\int_0^1\pi_1^*\om_1(\p_s v,J\p_s v)dtds\\
&\quad-\int_{-\infty}^{s_0}\int_0^1\frac{d}{ds}F_2(t,\pi_2\circ v)dtds-\int_{-\infty}^{s_0}\Big(\int_0^1\widetilde H_1(t,v)dt\Big)^2ds\\
&=-\int_{-\infty}^{s_0}\B(s)ds-\int_0^1F_2(t,\pi_2\circ v(s_0))-F_2(t,\pi_2\circ v_-)dt.\\
\eea
where $\B(s)$ is defined as
\bean
\B(s):=\int_0^1\pi_1^*\om_1(\p_s v,J\p_s v)dt+\Big(\int_0^1\widetilde H_1(t,v)dt\Big)^2.
\eea
Therefore the following estimate can be derived for any $s_0\in\R$
\bean
|\AA_1(w(s_0))|&\leq|\AA_1(w_+)|+2||F_2||_++\Big|\int_{-\infty}^{s_0}\B(s)ds\Big|,
\eea
and it remains to find a bound for $|\int_{-\infty}^{s_0}\B(s)ds|$.
Since $\B(s)$ is nonnegative, we are able to estimate as the following.
By setting $s_0=\infty$ in formula \eqref{eq:AA_1}, we have
\bean
\AA_1(w_+)-\AA_1(w_-)=-\int_{-\infty}^\infty\B(s)ds-\int_0^1F_2(t,\pi_2\circ v_+)-F_2(t,\pi_2\circ v_-)dt\\
\eea
Using the above formula, we obtain
\bean
\Big|\int_{-\infty}^{s_0}\B(s)ds\Big|&\leq\Big|\int_{-\infty}^{\infty}\B(s)ds\Big|\\
&\leq |\AA_1(w_+)|+|\AA_1(w_-)|+2||F_2||_+.
\eea
Thus we finally deduce
\beqn
|\AA_1(w(s_0))|\leq|\AA_1(w_+)|+2|\AA_1(w_-)|+4||F_2||_+,\qquad\forall s_0\in\R.
\eeq
\end{proof}

Once we have Lemma \ref{lemma:uniform bound of AA_i}, the rest of the proof of Theorem \ref{thm:compactness} is quite similar as in \cite{AF1}.

\begin{Lemma}\label{Lemma1}
Given a gradient flow line $w(s)=(v,\eta_1,\eta_2)(s)\in C^\infty(\R,\LLL\x\R^2)$ of $\AA^{\widetilde H_1,\widetilde H_2}_F$, assume that $v(t)\in U_\delta:=\widetilde G_1^{-1}(-\delta,\delta)\cap \widetilde G_2^{-1}(-\delta,\delta)$ for all $t\in(1/2,1)$ with $0<2\delta<\min\{1,\delta_0\}$. Then there exists $C_i>0$ satisfying
\beqn
|\eta_i|\leq C_i\Big(|\AA_i(v,\eta)|+||\nabla_m\AA^{\widetilde H_1,\widetilde H_2}_F||_m+1\Big),\qquad i=1,2.
\eeq
\end{Lemma}
\begin{proof}
We estimate
\bean
|\AA_i(v,\eta_1,\eta_2)|&=\Big|\int_0^1 v^*\pi_i^*\lambda_i+\eta_i\int_0^1\widetilde H_i(t,v)dt+\int_0^1F(t,v)dt\Big|\\
&\geq\Big|\eta_i\int_0^1\pi_i^*\lambda_i(v)\big(X_{\widetilde H_i}(t,v)\big)dt\Big|-\Big|\int_0^1\pi_i^*\lambda_i(v)\big(X_F(t,v)\big)dt\Big|-\Big|\eta_i\int^1_\frac{1}{2}\widetilde H_i(t,v)dt\Big|\\
&\,\,-\Big|\int_0^\frac{1}{2}F(t,v)dt\Big|-\Big|\int_0^1\pi_i^*\lambda_i(v)\big(\p_tv-\eta_1 X_{\widetilde H_1}(t,v)-\eta_2X_{\widetilde H_2}(t,v)-X_F(t,v)\big)dt\Big|\\
&\geq|\eta_i|-\delta|\eta_i|-C_{i,\delta}||\p_tv-\eta_1 X_{\widetilde H_1}(t,v)-\eta_2X_{\widetilde H_2}(t,v)-X_F(t,v)||_{L^1}-C_{i,\delta,F}\\
&\geq|\eta_i|-\delta|\eta_i|-C_{i,\delta}||\nabla_m\AA^{\widetilde H_1,\widetilde H_2}_F||_{m}-C_{i,F}\\
\eea
where $C_{i,\delta}:=||\pi_i^*\lambda_i|_{U_\delta}||_{L^\infty}$ and $C_{i,\delta,F}:=||F||_{L^\infty}+C_i||X_F||_{L^\infty}$. The second inequality holds since $\pi_i^*\lambda_i(X_{\widetilde H_j})=0$ if $i\ne j$. This estimate finishes the lemma with
\beqn
C_i:=\max\Big\{\frac{1}{1-\delta},\frac{C_{i,\delta}}{1-\delta},\frac{C_{i,\delta,F}}{1-\delta}\Big\},\qquad i=1,2.
\eeq
\end{proof}

\begin{Lemma}\label{Lemma2}
Given a gradient flow line $w(s)=(v,\eta_1,\eta_2)(s)\in C^\infty(\R,\LLL\x\R^2)$ of $\AA^{H_1,H_2}_F$, if there exists $t\in(\frac{1}{2},1)$ such that $v(t)\notin U_\delta$ then $||\nabla_m\AA^{\widetilde H_1,\widetilde H_2}_F(v,\eta_1,\eta_2)||_{m}>\epsilon$
for some $\epsilon=\epsilon_{\delta}$.
\end{Lemma}

\begin{proof}
Since $v(t)\notin U_\delta$ for some $t\in(\frac{1}{2},1)$, either $v(t)\notin U^1_\delta:=\widetilde G_1^{-1}(-\delta,\delta)$ or $v(t)\notin U^2_\delta:=\widetilde G_2^{-1}(-\delta,\delta)$ for that $t\in(\frac{1}{2},1)$. For simplicity, suppose $v(t)\notin U^1_\delta$.
If in addition $v(t)\notin U^1_{\delta/2}$ for all $t\in(\frac{1}{2},1)$, then we easily conclude that
\beqn
\bigr|\bigr|\nabla_m\AA^{\widetilde H_1,\widetilde H_2}_F(v,\eta_1,\eta_2)\bigr|\bigr|_{m}\geq\Big|\int_0^1\widetilde H_1(t,v(t))dt\Big|=\Big|\int^1_\frac{1}{2}\widetilde H_1(t,v(t))dt\Big|\geq\frac{\delta}{2}.
\eeq
Otherwise there is $t'\in(\frac{1}{2},1)$ such that $v(t')\in U^1_{\delta/2}$. Thus there exist $t_0,t_1\in(\frac{1}{2},1)$ satisfying one of the following two cases.
\bea\label{eq:first case}
v(t_0)\in \partial U^1_{\delta/2},\,\, v(t_1)\in \partial U^1_{\delta}
\quad\textrm{ and }\quad v(s)\in U^1_\delta-U^1_{\delta/2}\quad \textrm{   for all   } s\in [t_0,t_1]
\eea
or
\bean
v(t_1)\in \partial U^1_{\delta},\,\, v(t_0)\in \partial U^1_{\delta/2}
\quad\textrm{ and }\quad v(s)\in U^1_\delta-U^1_{\delta/2}\quad \textrm{   for all   } s\in [t_1,t_0].
\eea
We only treat the first case \eqref{eq:first case} and the second case follows analogously. With 
\beqn
\kappa:=\max_{x\in U_{\delta}}||\nabla_g \widetilde G_1(x)||_g
\eeq
we estimate
\bean
\kappa||\nabla_m\AA&^{\widetilde H_1,\widetilde H_2}_F(v,\eta_1,\eta_2)||_{m} \\
&\geq \kappa||\partial_t v-\eta_1 X_{\widetilde H_1}(t,v)-\eta_2X_{\widetilde H_2}(t,v)- X_{F}(t,v)||_{L^2}\\
&\geq \kappa||\partial_t v-\eta_1 X_{\widetilde H_1}(t,v)-\eta_2X_{\widetilde H_2}(t,v)-X_{F}(t,v)||_{L^1}\\
&\geq \int_{t_0}^{t_1}||\partial_t v-\eta_1 X_{\widetilde H_1}(t,v)-\eta_2X_{\widetilde H_2}(t,v)-X_{F}(t,v)||_g\cdot||\nabla \widetilde G_1(x)||_g dt\\
&\geq \bigg|\int_{t_0}^{t_1}\big\langle\nabla \widetilde G_1(v(t)),\partial_t v-\eta_1 X_{\widetilde H_1}(t,v)-\eta_2X_{\widetilde H_2}(t,v)-X_{F}(t,v)\big\rangle dt\bigg|\\
&= \bigg|\int_{t_0}^{t_1}d\widetilde G_1(v(t))\bigr(\partial_t v-\eta_1X_{\widetilde H_1}(t,v)-\eta_2X_{\widetilde H_2}(t,v)- \underbrace{X_{F}(t,v}_{=0}\bigr) dt\bigg|\\
&= \bigg|\int_{t_0}^{t_1}\frac{d}{dt}\widetilde G_1(v(t))dt-\underbrace{d \widetilde G_1(v(t))\big(\eta_1 X_{\widetilde H_1}(t,v)-\eta_2X_{\widetilde H_2}(t,v)\big)}_{=\eta_1\chi\om(X_{\widetilde G_1},X_{\widetilde G_1})+\eta_2\chi\om(X_{\widetilde G_1},X_{\widetilde G_2})=0}\bigg|\\
&\geq \big|\widetilde G_1(v(t_1))\big|-\big|\widetilde G_1(v(t_0))\big|\\
&=\frac{\delta}{2}.
\eea
Hence, the lemma follows with $\epsilon_\delta:=\min\bigr\{\delta/2,\delta/(2\kappa)\bigr\}.$\\
\end{proof}

Combining Lemma \ref{Lemma1} and  Lemma \ref{Lemma2}, we deduce the following {\em fundamental lemma}.
\begin{Lemma}\label{lemma:fundamental lemma}
For a gradient flow line $w(s)=(v,\eta_1,\eta_2)(s)\in C^\infty(\R,\LLL\x\R^2)$ of $\AA^{\widetilde H_1,\widetilde H_2}_F$, the following assertion holds for $i=1,2$ with some $C,\epsilon>0$.
\beqn
|\eta_i|\leq C\big(|\AA_i(w_-)|+|\AA_i(w_+)|+1\big) \quad\textrm{ if }\quad ||\nabla_m\AA^{\widetilde H_1,\widetilde H_2}_F(v,\eta_1,\eta_2)||_{m}<\epsilon.
\eeq
\end{Lemma}
\begin{proof}
According to Lemma \ref{Lemma2}, $v(t)$ lies in $U_\delta$ for all $t\in(\frac{1}{2},1)$ under the assumption
that $||\nabla_m\AA^{\widetilde H_1,\widetilde H_2}_F(v,\eta_1,\eta_2)||_m<\epsilon$. Thus we are able to apply Lemma \ref{Lemma1} and the following computation concludes the proof of the lemma.
\bean
|\eta_i|&\leq C_i(|\AA_i(v,\eta_1,\eta_2)|+||\nabla_m\AA^{\widetilde H_1,\widetilde H_2}_F(v,\eta_1,\eta_2)||_m+1)\\
&\leq C_i(2|\AA_i(w_-)|+|\AA_i(w_+)|+4||F_2||_++||\nabla_m\AA^{\widetilde H_1,\widetilde H_2}_F(v,\eta_1,\eta_2)||_m+1)\\
&\leq C_i(2|\AA_i(w_-)|+|\AA_i(w_+)|+4||F_2||_++1+\epsilon).
\eea
\end{proof}

\begin{Lemma}\label{lemma:tau(sigma)}
 For given a gradient flow line $w$ of $\AA^{\widetilde H_1,\widetilde H_2}_F$ and $\sigma\in\mathbb{R}$, we define
\beqn
\tau(\sigma):=\inf\bigr\{\tau\geq0\,\bigr|\,\,||\nabla_m\AA_{F}^{\widetilde H_1,\widetilde H_2}(w(\sigma+\tau))||_m\leq\epsilon\bigr\},
\eeq
Then we obtain a bound on $\tau(\sigma)$ as follows:
\beqn
\tau(\sigma)\leq\frac{\AA_{F}^{\widetilde H_1,\widetilde H_2}(w_-)-\AA_{F}^{\widetilde H_1,\widetilde H_2}(w_+)}{\epsilon^2}.
\eeq
\end{Lemma}
\begin{proof}
Using Lemma \ref{lemma:energy estimate for gradient lines}, we compute
\bean
\epsilon^2\tau(\sigma)
&\leq\int_\sigma^{\sigma+\tau(\sigma)}\big|\big|\nabla_m \AA_{F}^{\widetilde H_1,\widetilde H_2}(w)\big|\big|_m^2ds\\
&\leq E(w)\\
&\leq\AA^{\widetilde H_1,\widetilde H_2}_{F}(w_-)-\AA^{\widetilde H_1,\widetilde H_2}_{F}(w_+).
\eea
Dividing both sides through by $\epsilon^2$, the lemma follows.
\end{proof}

\begin{Thm}\label{thm:bound on eta}
Given two critical points $w_-$ and $w_+$, there exists a constant $\Theta>0$ depending only on $w_-$ and $w_+$
such that every gradient flow line $w(s)=(v,\eta_1,\eta_2)(s)$ of $\AA^{\widetilde H_1,\widetilde H_2}_F$ with fixed asymptotic ends $w_\pm$ satisfies
\beqn
||\eta_i||_{L^\infty}\leq\Theta \qquad\textrm{ for  } i=1,2.
\eeq
\end{Thm}
\begin{proof}
Using Lemma \ref{lemma:uniform bound of AA_i} and Lemma \ref{lemma:tau(sigma)}, we estimate
\bean
|\eta_i(\sigma)|&\leq|\eta_i(\sigma+\tau(\sigma))|+\int_\sigma^{\sigma+\tau(\sigma)}|\partial_s\eta_i(s)|ds\\
&\leq C\big(|\AA_i(w_-)|+|\AA_i(w_+)|+1\big)+\tau(\sigma)||\widetilde H_i||_{L^\infty}\\
&\leq C\big(|\AA_i(w_-)|+|\AA_i(w_+)|+1\big)+\Bigg(\frac{{\AA^{\widetilde H_1,\widetilde H_2}_{F}(w_-)-\AA^{\widetilde H_1,\widetilde H_2}_{F}(w_+)}}{\epsilon^2}\Bigg)||H_i||_{L^\infty}.
\eea
\end{proof}
As we mentioned before, Theorem \ref{thm:bound on eta} completes the proof of Theorem \ref{thm:compactness}.

\section{K\"unneth formula in Rabinowitz Floer homology}
Thanks to the previous section, we are now able to define Rabinowitz Floer homology of $(\Sigma_1\x\Sigma_2,M_1\x M_2)$ for admissible perturbations of the form $F_1\oplus F_2$ (or unperturbed). Whilst $\AA^{\widetilde H_1,\widetilde H_2}_F$ is generically Morse (Lemma \ref{lemma:genericity}), $\AA^{H_1,H_2}$ is never Morse because there is a $S^1$-symmetry coming from time-shift on the critical point set. However $\AA^{\widetilde H_1,\widetilde H_2}$ is generically Morse-Bott, so we are able to compute its Floer homology by choosing an auxiliary Morse function on the critical manifold and counting gradient flow lines with cascades, see \cite{F,CF}. Using the continuation method in Floer theory, we know that the Floer homology of $\AA^{\widetilde H_1,\widetilde H_2}_F$ is isomorphic to the  Floer homology of $\AA^{\widetilde H_1,\widetilde H_2}=\AA^{\widetilde H_1,\widetilde H_2}_0$ whenever these Floer homologies are defined. Thus we only treat the unperturbed Rabinowitz action functional $\AA^{\widetilde H_1,\widetilde H_2}$ and its Floer homology. Furthermore, we derive the K\"unneth formula by making use of the fact that all critical points and gradient flow lines can be split. In the last subsection, we prove Theorem B using similar steps to those in the proof of Theorem A; but we need to prove a special version of an isoperimetric inequality (Lemma \ref{lemma:isoperimetric ineq}) since unlike Theorem A, we have not insisted on any restrictions on perturbations in Theorem B.

\subsection{Rabinowitz Floer homology}
Firstly, we define a chain complex and a boundary operator for the Rabinowitz action functional. In order to define a chain complex we choose an additional Morse function $f$ on the critical manifold $\Crit\AA^{\widetilde H_1,\widetilde H_2}$. We define a $\Z/2$-Floer chain complex
\beqn
\CF_n(\AA^{\widetilde H_1,\widetilde H_2},f):=\Bigr\{\xi=\!\!\!\sum_{(v,\eta_1,\eta_2)}\!\!\xi_{(v,\eta_1,\eta_2)} (v,\eta_1,\eta_2)\,\Bigr|\,(v,\eta_1,\eta_2)\in\Crit_nf,\,\,\xi_{(v,\eta_1,\eta_2)}\in\Z/2\Bigr\}
\eeq
where $\xi_{(v,\eta_1,\eta_2)}$ satisfy the finiteness condition:
\beqn
\#\bigr\{{(v,\eta_1,\eta_2)}\in\Crit_nf\,\bigr|\,\xi_{(v,\eta_1,\eta_2)}\ne0,\,\,\AA^{H_1,H_2}(v,\eta_1,\eta_2)\geq\kappa\bigr\}<\infty,\quad\forall\kappa\in\R.
\eeq
The grading for the chain complex, $\mu=\mu_\RFH$, is described in the appendix of this paper.\\[-1.5ex]

To define the boundary operator, we roughly explain the notion of a {\em gradient flow line with cascades}. For rigorous and explicit constructions, we refer to \cite{F}. Consider a gradient flow line with cascades interchanging $w_-\subset C^-$ and $w_+\subset C^+$ where $w_\pm\in\Crit f$ and $C^\pm\subset\Crit\AA^{\widetilde H_1,\widetilde H_2}$; it starts with a gradient flow line of $f$ in $C^-$ with the negative asymptotic end $w_-$ and meets the negative asymptotic ends of a gradient flow line of $\AA^{\widetilde H_1,\widetilde H_2}$ (solving \eqref{eqn:gradient flow equation} with $F\equiv0$). We refer to this gradient flow line as a {\em cascade}. Its positive asymptotic end encounters a gradient flow line of $f$ in $C^+$ which converges to $w_+$. Several cascades and no cascades are also allowed. Now, we define a moduli space
\beqn
\widehat\M\{w_-,w_+\}:=\Bigg\{w\in C^\infty(\R,\LLL\x\R^2)\,\Bigg|\,\begin{array}{ll} w \textrm{ is a gradient flow line with cascades}\\[0.5ex]
\textrm{ with } \lim_{s\to\pm\infty}w(s)=w_\pm\in\Crit f\end{array}\Bigg\}
\eeq
and divide out the $\R$-action from shifting the gradient flow lines in the $s$-variable. Then we obtain the moduli space of unparametrized gradient flow lines, denoted by
\beqn
\M:=\widehat\M/\R.
\eeq

The standard transversality theory shows that this moduli space is a smooth manifold for a generic choice of the almost complex structure and the metric, see \cite{FHS,BO}. From the calculation \eqref{eq:dim, index computation} in the appendix, we also know that the dimension of $\M$ is equal to $\mu_\RFH(w_-)-\mu_\RFH(w_+)-1$. Therefore if $\mu_\RFH(w_-)-\mu_\RFH(w_+)=1$, $\M(w_-,w_+)$ is a finite set because of Theorem \ref{thm:compactness}. We let $\#_2\M\{w_-,w_+\}$ be the parity of  this moduli space. We define the boundary maps $\{\p^{\widetilde H_1,\widetilde H_2}_n\}_{n\in\Z}$ as follows:
\bean
\p^{\widetilde H_1,\widetilde H_2}_{n+1}:\CF_{n+1}(\AA^{\widetilde H_1,\widetilde H_2})&\pf\CF_{n}(\AA^{\widetilde H_1,\widetilde H_2})\\
w_-&\longmapsto\sum_{w_+\in\Crit_{n}f} \#_2\M\{w_-,w_+\}w_+.
\eea
Due to the Floer's central theorem, we know that $\p^{H_1,H_2}_n\circ\p^{H_1,H_2}_{n+1}=0$ so that $\big(\CF_*(\AA^{H_1,H_2}),\p^{H_1,H_2}_*\big)$ is a chain complex. We define Rabinowitz Floer homology by
\beqn
\RFH_n(\Sigma_1\x \Sigma_2,M_1\x M_2):=\HF_n\big(\AA^{\widetilde H_1,\widetilde H_2}\big)=\H_n\big(\CF_*(\AA^{\widetilde H_1,\widetilde H_2}),\p^{\widetilde H_1,\widetilde H_2}_*\big).
\eeq

\begin{Rmk}
Since in the previous section, we achieved the compactness result for $\AA^{\widetilde H_1,\widetilde H_2}_{F_1\oplus F_2}$, the Floer homology $\HF_n\big(\AA^{\widetilde H_1,\widetilde H_2}_{F_1\oplus F_2}\big)$ can be defined; besides, it is isomorphic to $\RFH_n(\Sigma_1\x \Sigma_2,M_1\x M_2)$ by the continuation homomorphism which counts gradient flow lines of $\AA^{\widetilde H_1,\widetilde H_2}_{F_s}$ where $F_s$ is a homotopy between $F_1\oplus F_2$ and $F\equiv0$.
\end{Rmk}
\subsection{Proof of Theorem A}
At first, we set
$$
H_1(t,x_1)=\chi(t)G_1(x_1)\in C^\infty(S^1\x M_1),\quad H_2(t,x_2)=\chi(t) G_2(x_2)\in C^\infty(S^1\x M_2)
$$
where $\chi:S^1\to [0,\infty)$ with $\int_0^1\chi dt=1$ and $\Supp\chi\subset(\frac{1}{2},1)$; it is clear that
\beqn
(\pi_i)_*X_{\widetilde H_i}(x_1,x_2)=X_{H_i}(x_i),\quad i=1,2.
\eeq
We consider the Rabinowitz action functionals $\AA^{H_1}:\LLL_{M_1}\x\R\to\R$ and $\AA^{H_2}:\LLL_{M_2}\x\R\to\R$:
\bean
&\bullet\quad\AA^{H_1}(v_1,\eta_1)=-\int_0^1 v_1^*\lambda_1-\eta_1\int_0^1H_1(t,v_1)dt,\\
&\bullet\quad\AA^{H_2}(v_2,\eta_2)=-\int_0^1 v_2^*\lambda_2-\eta_2\int_0^1H_2(t,v_2)dt.
\eea
In fact, we can accomplish compactness of gradient flow lines of each action functional with minor modifications of our case, or see \cite{AF1}. We observe that $(v_1,\eta_1)\in\Crit\AA^{H_1}$ solves
\beq\label{eq:crit1}
\partial_tv_1=\eta_1X_{H_1}(t,v_1)\quad\&\quad\int_0^1H_1(t,v_1)dt=0,
\eeq
and $(v_2,\eta_2)\in\Crit\AA^{H_2}$ solves
\beq\label{eq:crit2}
\partial_tv_2=\eta_2X_{H_2}(t,v_2)\quad\&\quad\int_0^1H_2(t,v_2)dt=0.
\eeq
Moreover a gradient flow line $w_1(s,t)=\big(v_1(s,t),\eta_1(s)\big):\R\x S^1\to M_1\x\R$ resp. $w_2(s,t)=\big(v_2(s,t),\eta_2(s)\big):\R\x S^1\to M_2\x\R$
is characterized by
\bea\label{eq:grad. flow eqs.}
\bullet\,\,&\p_sv_1+J_1(v_1)\big(\p_tv_1-\eta_1X_{H_1}(t,v_1)\big)=0,\quad\p_s\eta_1-\int_0^1H_1(t,v_1)dt=0,&&\quad\textrm{resp.} \\[1ex]
\bullet\,\,&\p_sv_2+J_2(v_2)\big(\p_tv_2-\eta_2X_{H_2}(t,v_2)\big)=0,\quad\p_s\eta_2-\int_0^1H_2(t,v_2)dt=0.&&
\eea
Then we define chain complexes $\CF(\AA^{H_1})$, $\CF(\AA^{H_2})$ and their boundary operators $\p^{H_1}$, $\p^{H_2}$ analogously as before, or see \cite{AF1} and denote their Floer homologies by
\beqn
\RFH(\Sigma_1,M_1)=\H\big(\CF(\AA^{H_1}),\p^{H_1}\big),\quad\RFH(\Sigma_2,M_2)=\H\big(\CF(\AA^{H_2}),\p^{H_2}\big).
\eeq
Next, for the K\"unneth formula, we define the tensor product of chain complexes by
\beqn
\bigr(\CF_*(\AA^{H_1})\otimes\CF_*(\AA^{H_2})\bigr)_n=\bigoplus_{i=0}^n\CF_i(\AA^{H_1})\otimes\CF_{n-i}(\AA^{H_2}).
\eeq
together with the boundary operator $\p_n^\otimes$ given by
\beqn
\p_n^\otimes\big((v_1,\eta_1)_i\otimes(v_2,\eta_2)_{n-i}\big)=\p_i^{H_1}(v_1,\eta_1)_i\otimes(v_2,\eta_2)_{n-i}+(v_1,\eta_1)_i\otimes\p_{n-i}^{H_2}(v_2,\eta_2)_{n-i}.
\eeq
Analyzing the critical point equations \eqref{eqn:critical point equation} when $F\equiv0$, \eqref{eq:crit1}, and \eqref{eq:crit2}, we easily notice that $\big((v_1,v_2),\eta_1,\eta_2\big)=(v,\eta_1,\eta_2)\in\Crit\AA^{H_1,H_2}$ if and only if $(v_1,\eta_1)\in\Crit\AA^{H_1}$ and $(v_2,\eta_2)\in\Crit\AA^{H_2}$ where $v_1=\pi_1\circ v:S^1\to M_1$ and $v_2=\pi_2\circ v:S^1\to M_2$ for the projections $\pi_1,\pi_2$.
Here, $(v_1,v_2)\in C^\infty(S^1,M_1\x M_2)$ is defined by
\bean
(v_1,v_2):S^1&\pf M_1\x M_2,\\
t&\longmapsto(v_1(t),v_2(t)).
\eea
Moreover the index behaves additively (see \eqref{eq:additicity of index}), thus we have
\beqn
\Crit_n(\AA^{\widetilde H_1,\widetilde H_2})=\bigcup_{i+j=n}\Crit_i(\AA^{H_1})\x\Crit_j(\AA^{H_2}),
\eeq
and we are able to define a chain homomorphism:
\bean
P_n:\bigr(\CF_*(\AA^{H_1})\otimes\CF_*(\AA^{H_2})\bigr)_n&\pf\CF_n(\AA^{\widetilde H_1,\widetilde H_2}),\\
(v_1,\eta_1)\otimes(v_2,\eta_2)&\longmapsto\big((v_1,v_2),\eta_1,\eta_2\big).
\eea
To verify that $P_n$ is a chain homomorphism, we need to show that
\beqn
\p_n^{H_1,H_2}\circ P_n=P_{n-1}\circ\p_n^\otimes.
\eeq
For $w_{1-}=(v_{1-},\eta_{1-})\in\Crit\AA^{H_1}$ and $w_{2-}=(v_{2-},\eta_{2-})\in\Crit\AA^{H_2}$, we compute
\bean
\p_n^{\widetilde H_1,\widetilde H_2}\circ P_n(w_{1-}\otimes w_{2-})&=\p_n^{\widetilde H_1,\widetilde H_2}\underbrace{\big((v_{1-},v_{2-}),\eta_{1-},\eta_{2-}\big)}_{=:w_-}\\
&=\!\!\!\!\!\sum_{\substack{w_+\in\Crit\AA^{\widetilde H_1,\widetilde H_2};\\\mu(w_+)=\mu(w_-)-1}}\!\!\!\!\!\#_2\M\{w_-,w_+\}w_+\\
&=\!\!\!\!\!\sum_{\substack{(v_{1+},\eta_{1+})\in\Crit\AA^{H_1};\\\mu(w_{1+})=\mu(w_{1-})-1}}\!\!\!\!\!\#_2\M\big\{w_-,((v_{1+},v_{2-}),\eta_{1+},\eta_{2-})\big\}\big((v_{1+},v_{2-}),\eta_{1+},\eta_{2-}\big)\\
&\,\,+\!\!\!\!\!\sum_{\substack{(v_{2+},\eta_{2+})\in\Crit\AA^{H_2};\\\mu(w_{2+})=\mu(w_{2-})-1}}\!\!\!\!\!\#_2\M\big\{w_-,((v_{1-},v_{2+}),\eta_{1-},\eta_{2+})\big\}\big((v_{1-},v_{2+}),\eta_{1-},\eta_{2+}\big)\\
&=\!\!\!\!\!\sum_{\substack{(v_{1+},\eta_{1+})\in\Crit\AA^{H_1};\\\mu(w_{1+})=\mu(w_{1-})-1}}\!\!\!\!\!\#_2\M\big\{w_{1-},w_{1+}\big\}P_{n-1}(w_{1+}\otimes w_{2-})\\
&\,\,+\!\!\!\!\!\sum_{\substack{(v_{2+},\eta_{2+})\in\Crit\AA^{H_2};\\\mu(w_{2+})=\mu(w_{2-})-1}}\!\!\!\!\!\#_2\M\big\{w_{2-},w_{2+}\big\}P_{n-1}(w_{1-}\otimes w_{2+})\\[0.5ex]
&=P_{n-1}(\p_i^{H_1}w_{1-}\otimes w_{2-})+P_{n-1}(w_{1-}\otimes\p_{n-i}^{H_2}w_{2-})\\[0.5ex]
&=P_{n-1}\circ\p_n^\otimes(w_{1-}\otimes w_{2-}).
\eea
where $\M\big\{w_{1-},w_{1+}\big\}$ resp. $\M\big\{w_{2-},w_{2+}\big\}$ is the moduli space which consists of gradient flow lines with cascades of $\AA^{H_1}$ resp. $\AA^{H_2}$. The fourth equality follows by comparing \eqref{eqn:gradient flow equation} together with \eqref{eq:grad. flow eqs.}.
Therefore we have an isomorphism
\beqn
(P_\bullet)_*:\H_\bullet\big(\CF(\AA^{H_1})\otimes\CF(\AA^{H_2})\big)\stackrel{\cong}\pf\H_\bullet(\CF(\AA^{\widetilde H_1,\widetilde H_2}))=\RFH_\bullet(\Sigma_1\x\Sigma_2,M_1\x M_2).
\eeq
Finally, the algebraic K\"unneth formula enable us to derive the desired (topological) K\"unneth formula in Rabinowitz Floer homology.
\beqn
\RFH_n(\Sigma_1\x\Sigma_2,M_1\x M_2)\cong\bigoplus_{p=0}^n\RFH_p(\Sigma_1,M_1)\otimes\RFH_{n-p}(\Sigma_2,M_2).\\[1ex]
\eeq
\subsection{Proof of Theorem B}
In this subsection, we do not consider $\Sigma_2$ and let $(M_2,\om_2)$ be closed and symplectically aspherical, i.e. $\om_2|_{\pi_2(M_2)}$. To prove Statement (B1) in Theorem B, we need compactness of gradient flow lines of the perturbed Rabinowitz action functional on $(\Sigma_1\x M_2,M_1\x M_2)$ for an arbitrary perturbation $F\in C_c^\infty(S^1\x M_1\x M_2)$. For that reason, we analyze the Rabinowitz action functional again as in Section 4; once we obtain the fundamental lemma, then the remaining steps are exactly same as before. Moreover due to compactness of gradient flow lines, we can find a leafwise intersection point for Hofer-small Hamiltonian diffeomorphisms using the stretching the neck argument in \cite{AF1}. We assume that $\Sigma_1\x M_2$ bounds a compact region in $M_1\x M_2$ for Statement (B2) in Theorem B throughout this subsection; but, when it comes to the existence of leafwise intersections, $\Sigma_1\x M_2$ need not bound a compact region in $M_1\x M_2$ using the techniques in \cite{Ka1,Ka2}. As before, we choose a defining Hamiltonian function $G\in C^\infty(M_1)$ so that
\begin{enumerate}
\item $G^{-1}(0)=\Sigma_1$ is a regular level set and $dG$ has a compact support.
\item $G_i(\phi_{Y}^t(x))=t$ for all $x\in\Sigma_i$, and $|t|<\delta$;
\end{enumerate}
where $Y$ is the Liouville vector field for $\Sigma_1\subset M_1$. We define $\widetilde G\in C^\infty(M_1\x M_2)$ by $\widetilde G(x_1,x_2)=G(x_1)$. Thus $\widetilde G$ is a defining Hamiltonian function for $\Sigma_1\x M_2$. We let $\widetilde H(t,x)=\chi(t)\widetilde G(x)\in C^\infty(S^1\x M_1\x M_2)$. With a perturbation $F\in C_c^\infty(S^1\x M_1\x M_2)$, the perturbed Rabinowitz action functional $\AA^{\widetilde H}_F:\LLL\x\R\pf\R$ is given by
\beqn
\AA^{\widetilde H}_F(v,\eta)=-\int_{D^2}\bar v^*\om_1\oplus\om_2-\eta\int_0^1\widetilde H(t,v)dt-\int_0^1F(t,v)dt
\eeq
where $\LLL=\LLL_{M_1\x M_2}\subset C^\infty(S^1,M_1\x M_2)$ is the component of contractible loops in $M_1\x M_2$ and $\bar v:D^2\to M_1\x M_2$ is a filling disk of $v$. The symplectic asphericity condition implies that the value of the above action functional is independent of the choice of filling disc.\\[-1ex]

Next, we prove the following lemma using a kind of isoperimetric inequality.
\begin{Lemma}\label{lemma:isoperimetric ineq}
Let $w(s,t)=(v(s,t),\eta(s))\in C^\infty(\R\x S^1,M_1\x M_2)\x C^\infty(\R,\R)$ be a gradient flow line of $\AA^{\widetilde H}_F$. We set $\gamma(t)=v(s_0,t)\in C^\infty(S^1,M_1\x M_2)$ for some fixed $s_0\in\R$. Then
$\int_{D^2}\bar \gamma^*\pi_2^*\om_2 $ is uniformly bounded provided $||\nabla_m\AA^{\widetilde H}_F(v(s_0,\cdot),\eta(s_0))||_m<\epsilon$ for some $\epsilon>0$:
\beq\label{eq:isoperimetric inequality}
\Big|\int_{D^2}\bar \gamma^*\pi_2^*\om_2\Big| \leq \max_{x\in\widetilde{M_2}}\bigr\{||\lambda_{\widetilde M_2}(x)||_{\tilde g_2}\,\bigr|\,d_{\tilde g_2}(x,\widetilde M_\star)<\epsilon+||X_F||_{L^\infty}\bigr\}\bigr(\epsilon+||X_F||_{L^\infty}\bigr).
\eeq
where $\widetilde{M_2}$ is the universal covering of $M_2$; $\tilde g_2$ is the lifting of the metric $g_2(\cdot,\cdot)=\om_2(\cdot,J_2\cdot)$ on $M_2$; $\widetilde M_\star$ is a fundamental domain in $\widetilde{M_2}$; $d_{\tilde g_2}(x,\widetilde M_\star)$ is the distance between $x$ and $\widetilde M_\star$; the value on the right hand side of \eqref {eq:isoperimetric inequality} is finite since $\widetilde M_\star\cong M_2$ is compact.
\end{Lemma}

\begin{proof}
We write $v(s,t)$ as $v(s,t)=(v_1,v_2)(s,t)$ where $v_1:\R\x S^1\to M_1$ and $v_2:\R\x S^1\to M_2$. Let $\gamma\in C^\infty(S^1,M_1\x M_2)$ be defined by $\gamma(t)=v(s_0,t)$ for some $s_0\in\R$. Since $\gamma$ is contractible and $M_2$ is symplectically aspherical, the value of $\int_{D^2}\bar\gamma^*\pi_2^*\om_2$ is well-defined. Let $\gamma_2:=\pi_2\circ\gamma$. We also consider $(\widetilde M_2,\widetilde{\om_2})$ the universal cover of $M_2$ where $\widetilde{\om_2}$ is the lift of $\om_2$ and we also lift the metric $g_2$ on $M_2$ which we write as $\tilde g_2$. Since we have assumed the symplectically asphericity of $(M_2,\om_2)$, there exists a primitive one form $\lambda_{\widetilde M_2}$ of $\widetilde{\om_2}$. Let $\widetilde M_\star(\cong M_2)$ be one of the fundamental domains in $\widetilde M_2$ and $\tilde v(s,t):\R\x S^1\to M_1\x\widetilde M_2$ be the lift of $v$ such that $\tilde v(s_0,t)=\tilde \gamma(t)$ intersects $M_1\x\widetilde M_\star$. Now, we can show the following kind of isoperimetric inequality. This inequality concludes the proof.
\bean
\Big|\int_{D^2}\bar \gamma^*\pi_2^*\om_2\Big| &=\Big|\int_{D^2}(\tilde{\bar \gamma}_2)^*\widetilde{\om_2}\Big|=\Big|\int_0^1\tilde \gamma_2^*\lambda_{\widetilde M_2}\Big|\\
&\leq||\lambda_{\widetilde M_2}|_{\gamma_2(S^1)}||_{L^\infty}\int_0^1||\p_t\tilde \gamma_2||_{\tilde g_2}dt\\
&=||\lambda_{\widetilde M_2}|_{\gamma_2(S^1)}||_{L^\infty}\int_0^1||\p_t \gamma_2||_{g_2}dt\\
&=||\lambda_{\widetilde M_2}|_{\gamma_2(S^1)}||_{L^\infty}\int_0^1||J\p_s \gamma_2+\pi_{2*}X_F(t,\gamma_2)||_{g_2}dt\\
&\leq\lambda_\mathrm{Max}\bigr(||\nabla_m\AA^{\widetilde H}_F(v(s_0,\cdot),\eta(s_0))||_m+||X_F||_{L^\infty}\bigr).
\eea
where
\bean
\lambda_\mathrm{Max}&:=\max_{x\in\widetilde{M_2}}\Big\{||\lambda_{\widetilde M_2}(x)||_{\tilde g_2}\,\Big|\,d_{\tilde g_2}(x,\widetilde M_\star)<\int_0^1||\p_t \gamma_2||_{g_2}dt\Big\}\\
&\leq\max_{x\in\widetilde{M_2}}\bigr\{||\lambda_{\widetilde M_2}(x)||_{\tilde g_2}\,\bigr|\,d_{\tilde g_2}(x,\widetilde M_\star)<||\nabla_m\AA^{\widetilde H}_F(v(s_0,\cdot),\eta(s_0))||_m+||X_F||_{L^\infty}\bigr\}.
\eea
\end{proof}

The following two lemmas can be proved similarly to the corresponding lemmas in the previous section.
\begin{Lemma}\label{lemma1}
We assume that for $(v,\eta)\in C^\infty(S^1,M_1\x M_2)\x\R$, $v(t)\in U_\delta:=\widetilde G^{-1}(-\delta,\delta)$ for all $t\in(\frac{1}{2},1)$ with $0<2\delta<\min\{1,\delta_0\}$. Then there exists $C>0$ satisfying
\beqn
|\eta|\leq C\Big(|\AA_F^{\widetilde H}(v,\eta)|+||\nabla_m\AA^{\widetilde H}_F(v,\eta)||_m+\Big|\int_{D^2}\bar v^*\pi^*_2\om_{2}\Big|+1\Big).
\eeq
\end{Lemma}

\begin{Lemma}\label{lemma2}
For $(v,\eta)\in C^\infty(S^1,M_1\x M_2)\x\R$ if there exists $t\in[\frac{1}{2},1]$ such that $v(t)\notin U_\delta$, then $||\nabla_m\AA^{\widetilde H}_F(v,\eta)||_{m}>\epsilon$
for some $\epsilon=\epsilon_{\delta}$.
\end{Lemma}

Due to the three previous lemmata, we deduce the fundamental lemma in the situation of Theorem B, and thus we obtain a uniform $L^\infty$-bound on the Lagrange multiplier $\eta$ by the same argument as in the previous section.
\begin{Lemma} For a gradient flow line $w(s)=(v,\eta)(s)\in C^\infty(\R,\LLL\x\R)$, the following assertions holds with some $C,\epsilon>0$. If $ ||\nabla_m\AA^{\widetilde H}_F(v,\eta)||_{m}<\epsilon$,
\beqn
|\eta|\leq C\big(|\AA_F^{\widetilde H}(w_-)|+|\AA_F^{\widetilde H}(w_+)|+\epsilon+\Xi_\epsilon+1\big) \quad\textrm{ provided that }\quad ||\nabla_m\AA^{\widetilde H}_F(v,\eta)||_{m}<\epsilon
\eeq
where $\Xi_\epsilon=\max\bigr\{||\lambda_{\widetilde {M_2}}(x)||_{\tilde g_2}\,|\,d_{{\tilde g_2}}(x,\widetilde M_\star)<\epsilon+||X_F||_{L^\infty}\bigr\}\bigr(\epsilon+||X_F||_{L^\infty}\bigr)<\infty.$
\end{Lemma}
\begin{proof}
The proof is almost same as the proof of Lemma \ref{lemma:fundamental lemma}. Since $||\nabla_m\AA^{\widetilde H}_F(v,\eta)||_{m}<\epsilon$, $v(t)\subset U_{\delta}$ for $t\in(\frac{1}{2},1)$ by Lemma \ref{lemma2}. Thus Lemma \ref{lemma:isoperimetric ineq} and Lemma \ref{lemma1} prove the lemma.
\end{proof}
This fundamental lemma proves compactness of gradient flow lines as before. Let $\phi\in\Ham_c(M_1\x M_2,\om_1\oplus\om_2)$ be a Hamiltonian diffeomorphism with the Hofer norm less than $\wp(\Sigma_1,\lambda_1)$. We  consider a moduli space of gradient flow lines of the Rabinowitz action functional perturbed by a special smooth family of Hamiltonian functions. Then in the boundary of this moduli space, there is a broken gradient flow line of which one asymptotic end gives rise to either a leafwise intersection point of $\phi$ or a closed Reeb orbit with period less than $||\phi||$. But since $||\phi||<\wp(\Sigma_1,\lambda_1)$, there is no such a closed Reeb orbit and hence we obtain a leafwise intersection point. This is so called  the {\em stretching the neck} argument, see \cite{AF1,Ka2}. Even further, there exists a leafwise intersection point even if $\Sigma_1\x M_2$ does not bound a compact region in $M_1\x M_2$ due to the arguments in \cite{Ka1,Ka2}. Next, we define the Rabinowitz Floer homology for $(\Sigma_1\x M_2,M_1\x M_2)$ in the same way as before and derive the K\"unneth formula in this situation. First of all, we consider another two action functionals $\AA^{H}:\LLL_{M_1}\x\R\to\R$ and $\AA:\LLL_{M_2}\to\R$ defined by
$$
\AA^{ H}(v_1,\eta):=-\int_0^1v_1^*\lambda_1-\eta\int_0^1 H(t,v)dt,\quad
\AA(v_2):=-\int_{D^2}\bar v_2^*\om_2.
$$
where $H(t,x)=\chi(t)G(x)\in C^\infty(S^1\x M_1)$. We note that $\AA^{\widetilde H}$ is defined on $\LLL_{M_1\x M_2}\x\R$.\\
As in the proof of Theorem A, we compare critical points of $\AA^{\widetilde H}$ and $\AA^{H}$ as follows.
\bean
\Crit_n(\AA^{\widetilde H})&=\bigcup_{i+j=n}\Crit_i(\AA^{H})\x\Crit_j(\AA).
\eea
Since $\Crit\AA$ consists of one component $M_2$, any gradient flow line with cascades of $\AA$ necessarily has zero cascades, and hence is simply a gradient flow line of an additional Morse function $f\in C^\infty(M_2)$. Thus the chain group for the Morse-Bott homology of $\AA$ is given by $\CF(\AA,f)=\CM(f)$. Here $\CM$ stands for the Morse complex. The following map is a chain isomorphism, which can be verified using the methods of the previous subsection.
\bean
P_n:\bigr(\CF_*(\AA^{H})\otimes\CM_*(f)\bigr)_n&\pf\CF_n(\AA^{\widetilde H}),\\
(v_1,\eta)\otimes v_2&\longmapsto\big((v_1,v_2),\eta\big).
\eea
Therefore it induces an isomorphism on the homology level:
\beqn
(P_\bullet)_*:\H_\bullet\big(\CF(\AA^{H})\otimes\CM(f)\big)\stackrel{\cong}\pf\H_\bullet\big(\CF(\AA^{\widetilde H})\big)=\RFH_\bullet(\Sigma_1\x M_2,M_1\x M_2).
\eeq
Finally, the K\"unneth formula for $(\Sigma_1\x M_2,M_1\x M_2)$ directly follows:
\beqn
\RFH_n(\Sigma_1\x M_2,M_1\x M_2)\cong\bigoplus_{p=0}^n\RFH_p(\Sigma_1,M_1)\otimes\H_{n-p}(M_2).\\[1ex]
\eeq

\section{Applications}
As we have mentioned in the introduction, we cannot achieve compactness of gradient flow lines of $\AA^{\widetilde H_1,\widetilde H_2}_F$ for an arbitrary perturbation $F$. For that reason, the existence problem of leafwise intersection points for a product submanifold which is not of contact type is still open. On the other hand, the existence of leafwise intersection points for contact coisotropic submanifolds was already proved in \cite{Gu,Ka2}. Furthermore, due to the K\"unneth formula, we can deduce the existence of infinitely many leafwise intersection points for some kind of product submanifolds of contact type.  First, we recall the notion of contact condition on coisotropic submanifolds introduced by Bolle \cite{Bo1,Bo2}.

\begin{Def}
A coisotropic submanifold $\Sigma$ of codimension $k$ in a symplectic manifold $(M,\om)$ is called of {\em restricted contact type} if there exist global one forms $\lambda_1,\dots,\lambda_k\in\Omega^1(M)$ which satisfy
\begin{enumerate}
\item $d\lambda_i=\om$ for $i=1,\dots,k$;
\item $\lambda_1\wedge\cdots\wedge\lambda_k\wedge\om^{n-k}|_\Sigma\ne0$.
\end{enumerate}
\end{Def}
\begin{Rmk}\label{rmk:restriction of contact}\cite{Bo2,Gi1}
Let $\Sigma$ be closed and have contact type in $M$. Then a one form $\lambda=a_1\lambda_1+\cdots+a_k\lambda_k$ with $a_1+\cdots+a_k=0$ is closed and hence defines an element of $\H^1_\mathrm{dR}(\Sigma)$. In addition, $\lambda\neq0$ is not exact; otherwise $\lambda=df$ for some $f\in C^\infty(\Sigma)$, and hence $\lambda(x)=0$ at a critical point $x$ of $f$, but condition (ii) yields that $\lambda_1,\dots,\lambda_k$ are linearly independent on $\Sigma$; thus $\lambda_1(x)=\cdots\lambda_k(x)=0$. As a result, $\dim\H^1_\mathrm{dR}(\Sigma)\geq k-1$. It imposes restrictions on the contact condition; for instance, $S^3\x S^3$ is not of contact type in $\R^8$.
\end{Rmk}
We note that if the codimension of $\Sigma$ is bigger than one, $\Sigma$ never bounds a compact region in $M$. In spite of such a dimension problem, the condition that global coordinates exist (roughly speaking, Poisson-commuting Hamiltonian functions whose common zero locus is only $\Sigma$) enable us to unfold the generalized Rabinowitz Floer homology theory \cite{Ka2}. It turns out that a product of contact hypersurfaces bounding respective ambient symplectic manifolds has global coordinates.
\begin{Thm}\cite{Ka2}
If $\Sigma$ is a contact coisotropic submanifold of $M$ which admits global  coordinates, then the Floer homology of the perturbed Rabinowitz action functional is well-defined.
\end{Thm}
Since the Rabinowitz action functional can be defined for each homotopy classes of loops, we can define the Rabinowitz Floer homology $\RFH(\Sigma,M,\gamma)$ for $\gamma\in [S^1,M]$. We note that the  $\RFH(\Sigma,M)$ considered so far, is equal to $\RFH(\Sigma,M,\mathrm{pt})$. Moreover we also define Rabinowitz Floer homology on the full loop space and denote it by $\mathbf{RFH}(\Sigma,M)$. Then we have
\beqn
\mathbf{RFH}_*(\Sigma,M)=\bigoplus_{\gamma\in[S^1,M]}\RFH_*(\Sigma,M,\gamma).
\eeq
We recall the computation of Rabinowitz Floer homology on the (unit) cotangent bundle $(S^*N,T^*N)$ for a closed Riemannian manifold $N$.
\begin{Thm}\cite{CFO,AS,Me}
\beqn
\mathbf{RFH}_*(S^*N,T^*N) \cong \left\{ \begin{array}{ll}
 \H_*(\Lambda N), & *>1,\\[1ex]
 \H^{-*+1}(\Lambda N), & *<0. \end{array}\right.\;
\eeq
Here $\Lambda N$ stands for the free loop space of $N$.
\end{Thm}

Since the K\"unneth formula obviously holds for $\mathbf{RFH}$ also, the following corollary directly follows.
\begin{Cor}
If $\mathbf{RFH}_*(\Sigma_1,M_1)\neq 0$, and $\dim\H_*(\Lambda N)=\infty$ then
\beqn
\dim\mathbf{RFH}_*(\Sigma_1\x S^*N,M_1\x T^*N)=\infty.
\eeq
Accordingly, if $\Sigma_1\x S^*N$ has contact type again, $\Sigma_1\x S^*N$ has infinitely many leafwise intersection points or a periodic leafwise intersections for a generic perturbation.
\end{Cor}

\subsection{Proof of Corollary A and B}

From now on, we investigate leafwise intersections on $(S^*S^1\x S^*N,T^*S^1\x T^*N)$.
\begin{Lemma}\label{lemma:SS1 SN is contact}
$S^*S^1\x S^*N$ is a contact submanifold of codimension two in $T^*S^1\x T^*N$.
\end{Lemma}
\begin{proof}
$(T^*S^1,\om_{\mathrm{S^1,can}})\cong(S^1\x\R,d\theta\wedge dr)$ where $\theta$ is the angular coordinate on $S^1$ and $r$ is the coordinate on $\R$. Then $d\theta\wedge dr$ has two global primitives $-rd\theta$ and $-rd\theta+d\theta$. We can easily check that $S^*S^1\x S^*N$ carries a  contact structure with $-rd\theta\oplus\lambda_{\mathrm{N,can}}$ and $(-rd\theta+d\theta)\oplus\lambda_{\mathrm{N,can}}$ where $\lambda_{\mathrm{N,can}}$ is the canonical one form on $T^*N$.
\end{proof}
To exclude periodic leafwise intersections, we consider the loop space $\Omega$ defined by
\beqn
\Omega:=\bigr\{v=(v_1,v_2)\in C^\infty(S^1,T^*S^1\x T^*N)\,\bigr|\,v_1\textrm{ is contractible in } T^*S^1\bigr\}.
\eeq
Then we define the Rabinowitz action functional on this loop space, $\AA^{H_1,H_2}_F:\Omega\x\R^2\pf\R$, and construct the respective Rabinowitz Floer homology $\RFH(S^*S^1\x S^*N,T^*S^1\x T^*N,\Omega)$ as before. Moreover the following type of the K\"unneth formula holds.
\beqn
\RFH_n(S^*S^1\x S^*N,T^*S^1\x T^*N,\Omega)\cong\bigoplus_{p=0}^n\RFH_p(S^*S^1,T^*S^1)\otimes\mathbf{RFH}_{n-p}(S^*N,T^*N).
\eeq
Therefore $\RFH(S^*S^1\x S^*N,T^*S^1\x T^*N,\Omega)$ is of infinite dimensional if $\dim\H_*(\Lambda N)=\infty$ and Lemma \ref{lemma:genericity} below yields that there are infinitely many leafwise intersection points for a generic perturbation whenever $\dim N\geq2$. This proves Corollary A.

In order to prove that there is generically no periodic leafwise intersections, we review the argument in \cite{AF2}. We denote by $\mathcal{R}$ the set of closed Reeb orbits in $T^*N$ which has dimension one. It is convenient to introduce the following sets:
\beqn
\FF^j=\bigr\{F\in C^j_c(S^1\x T^*S^1\x T^*N)\,\bigr|\, F(t,\cdot)=0,\,\,\forall t\in\big[\frac{1}{2},1\big]\bigr\},\quad
\FF=\bigcap_{j=1}^\infty\FF^j.
\eeq

\begin{Lemma}\label{lemma:genericity}
If $\dim N\geq 2$, then the set
\beqn
\FF_{S^*S^1\x S^*N}:=\left\{F\in\FF\,\,\Bigg|\,\begin{aligned}\AA^{H_1,H_2}_F \textrm{ is Morse }\&\,\, v(0)\cap (S^*S^1\x R)=\emptyset\\\textrm{for all   } \forall(v,\eta_1,\eta_2)\in\AA^{H_1,H_2}_F,\,\,  R\in\mathcal{R}\end{aligned}\right\}\;
\eeq
is dense in the set $\FF$.
\end{Lemma}
\begin{proof}

In this proof, we denote by
\beqn
\Omega^{1,2}:=\bigr\{v=(v_1,v_2)\in W^{1,2}(S^1,T^*S^1\x T^*N)\,\bigr|\,v_1\textrm{ is contractible in } T^*S^1\bigr\}.
\eeq
the loop space which is indeed a Hilbert manifold. Let $\EE$ be the $L^2$-bundle over $\Omega^{1,2}$ with $\EE_v=L^2(S^1,v^*T(S^*S^1\x S^*N))$. We consider the section
\beqn
S:\Omega^{1,2}\x\R^2\x\FF^j\pf\EE^\vee\x\R^2\quad\textrm{defined by}\quad S(v,\eta_1,\eta_2,F):=d\AA^{H_1,H_2}_F(v,\eta_1,\eta_2).
\eeq
Here the symbol $\vee$ represents the dual space. At $(v,\eta_1,\eta_2,F)\in S^{-1}(0)$, the vertical differential
\beqn
DS:T_{(v,\eta_1,\eta_2,F)}\Omega^{1,2}\x\R^2\x\FF^j\pf\EE_v^\vee\x\R^2
\eeq
is given by the pairing
\beqn
\big\langle DS_{(v,\eta_1,\eta_2,F)}[\hat v^1,\hat\eta^1_1,\hat\eta_2^1,\hat F],[\hat v^2,\hat\eta^1_2,\hat\eta^2_2]\big\rangle=\mathscr{H}_{\AA^{H_1,H_2}_F}[(\hat v^1,\hat\eta^1_1,\hat\eta_2^1),(\hat v^2,\hat\eta^2_1,\hat\eta_2^2)]+\int_0^1\hat F(t,v)dt.
\eeq
where $\mathscr{H}_{\AA^{H_1,H_2}_F}$ is the Hessian of $\AA^{H_1,H_2}_F$. Due to the arguments in \cite{AF1} (in fact they proved the surjectivity for $\AA_F^H$, but their proof obviously can be extended to our situation, see also \cite{Ka2}), we know that for $(v,\eta_1,\eta_2,F)\in S^{-1}(0)$, $DS_{(v,\eta_1,\eta_2,F)}$ is surjective on the space
\beqn
\mathcal{V}:=\big\{(\hat v,\hat\eta_1,\hat\eta_2,\hat F)\in T_{(v,\eta_1,\eta_2,F)}(\Omega^{1,2}\x\R^2\x\FF^j)\,\big|\,\hat v(0)=0\big\}.
\eeq
Next, we consider the evaluation map
\bean
\mathrm{ev}:\M&\pf S^*S^1\x S^*N,\\
(v,\eta_1,\eta_2,F)&\longmapsto v(0).
\eea
Since $DS_{(v,\eta_1,\eta_2,F)}|_\mathcal{V}$ is surjective, Lemma \ref{lemma:Dietmar} below implies that $\mathrm{ev}$ is a submersion. Then $\MM_\mathcal{R}:=\mathrm{ev}^{-1}(S^*S^1\x\mathcal{R})$ is a submanifold in $\MM$ of
\beqn
\codim(\MM_\mathcal{R}/\MM)=\codim(S^*S^1\x\mathcal{R}/S^*S^1\x S^*N).
\eeq
We consider the projections $\Pi:\MM\pf\FF^j$ and $\Pi_\mathcal{R}:=\Pi_{|\MM_\mathcal{R}}$. Then $\AA^{H_1,H_2}_F$ is Morse if and only if $F$ is a regular value of $\Pi$, which is a generic property by Sard-Smale theorem (for $j$ large enough). The set $\Pi^{-1}(F)$ of leafwise intersection points for $F$ is manifold of required dimension zero since it is a critical set of $\AA^{H_1,H_2}_F$. On the other hand, $\Pi^{-1}_\mathcal{R}(F)$ is a manifold of dimension
\beqn
0+\dim\MM_\mathcal{R}-\dim\MM=-\codim(\MM_\mathcal{R}/\MM)<0
\eeq
since we have assumed $\dim N\geq 2$. Therefore $\mathrm{ev}$ does not intersect $S^*S^1\x\mathcal{R}$, so the set
\beqn
\FF_{S^*S^1\x S^*N}^j:=\FF_{S^*S^1\x S^*N}\cap\FF^j
\eeq
is dense in $\FF$ for all $j\in\N$. Since $\FF_{S^*S^1\x S^*N}$ is the  countable intersection of $\FF_{S^*S^1\x S^*N}^j$ for $j\in\N$, it is dense again in $\FF$ and the lemma is proved.
\end{proof}

\begin{Lemma}\label{lemma:Dietmar}\textbf{(Salamon)}
Let $\EE\pf\BB$ be a Banach bundle and $s:\BB\pf\EE$ a smooth section. Moreover, let $\phi:\BB\pf N$ be a smooth map into the Banach manifold $N$.
We fix a point $x\in s^{-1}(0)\subset\BB$ and set $K:=\ker d\phi(x)\subset T_x\BB$ and assume the following two conditions.
\begin{enumerate}
\item The vertical differential $Ds|_K:K\pf\EE_x$ is surjective.
\item $d\phi(x):T_x\BB\pf T_{\phi(x)}N$ is surjective.
\end{enumerate}
Then  $d\phi(x)|_{\ker Ds(x)}:\ker Ds(x)\pf T_{\phi(x)}N$ is surjective.
\end{Lemma}
\begin{proof}
Given $\xi\in T_{\phi(x)}N$, condition (ii) implies that there exists $\eta\in T_x\BB$ satisfying $d\phi(x)\eta=\xi$. In addition, by condition (i),
there exists $\zeta\in K\subset T_x\BB$ satisfying $Ds(x)\zeta=Ds(x)\eta$. We set $\tau:=\eta-\zeta$ and compute
\beqn
Ds(x)\tau=Ds(x)\eta-Ds(x)\zeta=0
\eeq
thus, $\tau\in\ker Ds(x)$. Moreover,
\beqn
d\phi(x)\tau=d\phi(x)\eta-\underbrace{d\phi(x)\zeta}_{=0}=d\phi(x)\eta=\xi
\eeq
proves the lemma.
\end{proof}

In the case of Theorem B, we also redefine the Rabinowitz action functional $\AA^H_F:\Omega_{M_2}\x\R\to\R$ by
\beqn
\AA^H_F(v,\eta)=-\int_0^1 v_1^*\lambda_1-\int_{D^2}\bar v_2^*\om_2-\eta\int_0^1H(t,v)dt-\int_0^1F(t,v)dt
\eeq
where
\beqn
\Omega_{M_2}:\big\{v=(v_1,v_2)\in C^\infty(S^1,M_1\x M_2)\,\big|\,v_2 \textrm{ is contractible in }M_2\big\}.
\eeq
We can also define the respective Rabinowitz Floer homology and derive an appropriate K\"unneth formula as before.
\begin{Cor}
Let $(M_2,\om_2)$ be a closed, symplectically aspherical symplectic manifold. If a closed manifold $N$ has $\dim\H_*(\Lambda N)=\infty$, then
\beqn
\dim\RFH_*(S^*N\x M_2,T^*N\x M_2,\Omega_{M_2})=\infty.
\eeq
Therefore, if $\dim N\geq2$, $S^*N\x M_2$ has infinitely many leafwise intersection points for a generic perturbation.
\end{Cor}
The previous corollary proves Corollary B. 

\begin{Rmk}
The corollaries still holds when we deal with a generic fiber-wise star shaped hypersurface $\Sigma\subset T^*N$ instead of $S^*N$, see \cite{AF2}.
\end{Rmk}

\section{Appendix : Index for Rabinowitz Floer homology}

In fact, we are able to derive the K\"unneth formula and obtain applications without defining indices. Nevertheless, for the sake of completeness, we briefly recall the index for generators of the Rabinowitz Floer chain complex in this appendix (see \cite{CF} for the detailed arguments). Let $\Sigma$ be a contact hypersurface in $M$. Under the following assumption the Rabinowitz Floer homology has $\Z$-grading,
\begin{enumerate}
\item[(H1)] Closed Reeb orbits on $(\Sigma,\lambda)$ is of Morse-Bott type \cite{CF}.\\[-0.8ex]
\item[(H2)] The first chern class $c_1$ vanishes on $TM$.
\end{enumerate}

\begin{Rmk}
Without any hypothesis on the first chern class, the Rabinowitz Floer homology has $\Z/2$-grading. The non-degeneracy assumption (H1) is satisfied for generic hypersurfaces and the invariance property allows us to perturb a hypersurface to be Morse-Bott type.
\end{Rmk}

Let $\M$ be the moduli space of all finite energy gradient flow lines of $\AA^H$ and $w=(v,\eta)\in C^\infty(\R\x S^1,M)\x C^\infty(\R,\R)$ be a gradient flow line of $\AA^H$ with $\lim_{s\to\pm\infty}w(s)=w_\pm=(v_\pm,\eta_\pm)\in \Crit f$ and $v_\pm\subset C_\pm$ where $C_\pm\subset\Crit\AA^H$ are connected components of the critical manifold and $f$ is an additional Morse function on a critical manifold $\Crit\AA^H$. The linearization of the gradient flow equation along $(v,\eta)$ gives rise to an operator $D_{(v,\eta)}^{\AA^H}$. For suitable weighted Sobolev spaces, $D_{(v,\eta)}^{\AA^H}$ is a Fredholm operator. Then the local virtual dimension of $\M$ at $(v,\eta)$ is defined to be
\beqn
\mathrm{virdim}_{(v,\eta)}\M:=\ind D_{(v,\eta)}^{\AA^H}+\dim C_-+\dim C_+.
\eeq
Here $\ind D_{(v,\eta)}^{\AA^H}$ stands for the Fredholm index of $D_{(v,\eta)}^{\AA^H}$. Cieliebak-Frauenfelder \cite{CF} investigated the spectral flow of the Hessian $\mathrm{Hess}_{\AA^H}$ and consequently proved the following index formula:
\bean
\mathrm{virdim}_{(v,\eta)}\M =\mu_{\mathrm{CZ}}(v_+)-\mu_{\mathrm{CZ}}(v_-)+\frac{\dim C_-+\dim C_+}{2}.
\eea
Here $\mu_{\mathrm{CZ}}$ is the {\em Conley-Zehnder index} defined below. Since a closed Reeb orbit $v_+$ is contractible in $M$, we have a filling disk $\bar v_+:D^2\to M$ such that $\bar v_+|_{\p D^2}=v_+$. The filling disk $\bar v_+$ determines homotopy class of trivialization of the symplectic vector bundle $(\bar v_+)^* TM$. The linearized flow of the Reeb vector field along $v_+$ defines a path in $Sp(2n,\R)$ the group of symplectic matrices. The Conley Zehnder index of $v_+$ is defined by the Maslov index of \cite{RS} this path. This index is independent of choice of filling disk due to (H2) because the Conley Zehnder indices of different filling disks are differ by $c_1$. In the same way $\mu_{\mathrm{CZ}}(v_-)$ is also defined. Now, we are in a position to define a grading $\mu_{\mathrm{RFH}}$ on $\CF(\AA^H)$ by
\beqn
\mu_{\mathrm{RFH}}(v_\pm,\eta_\pm):=\mu_{\mathrm{CZ}}(v_\pm)+\mu^f_\sigma(v_\pm).
\eeq
where $\mu^f_\sigma$ is the {\em signature index} defined by
\bean
\mu_f^\sigma(v_\pm)&=-\frac{1}{2}\mathrm{sign}(\mathrm{Hess}_f(v_\pm))\\
&=-\frac{1}{2}\Bigg(\#\left\{\begin{array}{l}\textrm{positive eigenvalues of}\\ \textrm{the Hessian of $f$ at $v_\pm$}\end{array}\right\}\;-\#\left\{\begin{array}{l}\textrm{negative eigenvalues of}\\ \textrm{the Hessian of $f$ at $v_\pm$}\end{array}\right\}\Bigg).
\eea
We note that by definition,
\beqn
\mu_f^\sigma(v_\pm)=\mu_f^{\mathrm{Morse}}(v_\pm)-\frac{1}{2}\dim C_\pm.
\eeq

Then we notice that the dimension of gradient flow lines of $\AA^H_F$ interchanging $w_-$ and $w_+$ equals the index difference of the two critical points by the following computation.
\bea\label{eq:dim, index computation}
\mathrm{dim}\widehat\M\{w_-,w_+\}&=\mathrm{virdim}_{(v,\eta)}\M-\dim C_+-\dim C_-+\dim W_f^u(v_-)+\dim W_f^s(v_+)\\
&=\mu_{\mathrm{CZ}}(v_-)-\mu_{\mathrm{CZ}}(v_+)-\frac{\dim C_-+\dim C_+}{2}+\mu_f^{\mathrm{Morse}}(w_-)+\dim C^+-\mu_f^{\mathrm{Morse}}(w_+)\\
&=\mu_{\mathrm{CZ}}(v_-)-\mu_{\mathrm{CZ}}(v_+)-\frac{\dim C_-+\dim C_+}{2}+\mu_f^\sigma(v_-)+\frac{1}{2}\dim C_-+\dim C_+\\
&\quad-(\mu^f_\sigma(v_+)+\frac{1}{2}\dim C_+)\\
&=\mu_{\mathrm{CZ}}(v_-)-\mu_{\mathrm{CZ}}(v_+)+\mu^f_\sigma(v_-)-\mu_\sigma(v_+)\\
&=\mu_{\mathrm{RFH}}(v_-)-\mu_{\mathrm{RFH}}(v_+)
\eea
where  $W_f^s(v^+)$($W_f^u(v^-)$) is the (un)stable manifold with respect to $(f,v^\pm)$. \\

Furthermore, the $\RFH$-index of $\big((v_1,v_2),\eta_1,\eta_2\big)\in\Crit\AA^{\widetilde H_1,\widetilde H_2}$ (as used in Theorem A) splits into the indices of $(v_1,\eta_1)$ and $(v_2,\eta_2)$.
\beq\label{eq:additicity of index}
\mu_{\mathrm{RFH}}\big((v_1,v_2),\eta_1,\eta_2\big)=\mu_{\mathrm{RFH}}(v_1,\eta_1)+\mu_{\mathrm{RFH}}(v_2,\eta_2)
\eeq
since the Conley-Zehnder index (in fact, the Maslov index) and the Morse index behave additively under the direct sum operation.

\subsection*{Acknowledgments} I am deeply indebted to Urs Frauenfelder for numerous fruitful discussions. I also thanks the anonymous referee for careful reading and comments.

\end{document}